\documentclass[11pt]{article}%
\usepackage{amsmath}
\usepackage{amsfonts}
\usepackage{amssymb}
\usepackage{graphicx}%
\newtheorem{theorem}{Theorem}

\newtheorem{corollary}[theorem]{Corollary}

\newtheorem{proposition}[theorem]{Proposition}

\begin{document}

\begin{center}
{\Large On Non-occurrence of Chaos in Nonautonomous Planar Flows}

\medskip

\centerline{H. SEDAGHAT \footnote{Department of Mathematics, Virginia Commonwealth University Richmond, Virginia, 23284-2014, USA; Email: h.sedagha@vcu.edu}}
\end{center}

\medskip

\begin{abstract}
%\noindent\textbf{Abstract}
By folding nonautonomous differential systems in the plane to scalar differential
equations, a sufficient condition for the non-occurrence of chaotic behavior
is obtained.

\end{abstract}

\bigskip

It is a well-known fact that autonomous (time-independent) planar flows do not
exhibit the type of complex aperiodic behavior known as chaos \cite{P}.
It is also widely known that nonautonomous (time-dependent) flows may be
chaotic, as in the case of the periodically forced Duffing's equation
\cite{ST}. Being nonautonomous is a necessary but not a sufficient condition
for the occurrence of chaos in the plane; for instance, a nonhomogeneous
linear system with constant coefficients has no chaotic flows. Nonlinear,
time-dependent planar systems may also fail to be chaotic although sufficient
conditions for the occurrence or non-occurrence of chaos in the plane are not
readily available for nonlinear systems.

In this paper, we obtain a sufficient condition for the non-occurrence of
chaotic planar flows by folding the system into an ordinary differential
equation of order 2. If the differential equation does not have chaotic
solutions then neither does the original system. Folding has been used in
control theory (the controllability canonical form--see, e.g., \cite{B}) and
in deriving jerk functions for autonomous differential systems of three
equations; see \cite{E}. It is defined as a general algorithmic process in
\cite{SF}.

Consider the differential system
\begin{equation}
\left\{
\begin{array}
[c]{c}%
\dot{x}=f(t,x(t),y(t))\\
\dot{y}=g(t,x(t),y(t))
\end{array}
\right.  \label{ds1}%
\end{equation}
where $\dot{x}=dx/dt,\dot{y}=dy/dt$ and each of the functions $f,g:I\times
D\rightarrow(-\infty,\infty)$ has continuous partial derivatives with respect
to $t,x,y$. Here $D$ is an invariant domain in the plane $\mathbb{R}^{2}$ and
$I$ is a given interval.

Suppose that $f$ satisfies the condition
\begin{equation}
\dot{x}(t)=f(t,x(t),y(t))\Rightarrow y(t)=h(t,x(t),\dot{x}(t))\quad\text{for
all }t\label{des}%
\end{equation}
i.e., $f$ is semi-invertible with respect to $y$; see \cite{SF}. For instance,
separable functions such as%
\[
f(t,x(t),y(t))=f_{1}(t,x(t))+f_{2}(t,y(t))\quad\text{or\quad}%
f(t,x(t),y(t))=f_{1}(t,x(t))f_{2}(t,y(t))
\]
are semi-invertible with respect to $y$ if and only if $f_{2}$ is injective
and there is $f_{2}^{-1}$ such that $f_{2}^{-1}(t,f_{2}(t,y(t)))=y(t).$ In
particular, all affine functions $f(t,x(t),y(t))=a(t)x(t)+b(t)y(t)+c(t)$ are
separable; further, they are semi-invertible with respect to $y$ if and only
if $b(t)\not =0$ for all $t.$

By the chain rule,%
\begin{align*}
\overset{..}{x} &  =\frac{d}{dt}f(t,x,y)=f_{t}(t,x,y)+\dot{x}f_{x}%
(t,x,y)+\dot{y}f_{y}(t,x,y)\\
&  =f_{t}(t,x,y)+\dot{x}f_{x}(t,x,y)+g(t,x,y)f_{y}(t,x,y)
\end{align*}

Now using (\ref{des}),%
\begin{equation}
\overset{..}{x}=f_{t}(t,x,h(t,x,\dot{x}))+\dot{x}f_{x}(t,x,h(t,x,\dot
{x}))+g(t,x,h(t,x,\dot{x}))f_{y}(t,x,h(t,x,\dot{x})). \label{de2}%
\end{equation}

Let $\phi(t,x,\dot{x})$ denote the right hand side of (\ref{de2}). For an
initial point $(x(t_{0}),y(t_{0}))\in D$ where $t_{0}\in I$, the ordinary
differential equation%
\begin{equation}
\overset{..}{s}=\phi(t,s,\dot{s}) \label{de1}%
\end{equation}
with initial values $s(t_{0})=x(t_{0}),\ \dot{s}(t_{0})=f(t_{0},x(t_{0}%
),y(t_{0}))$ gives the x-component of the flow $(x(t),y(t))$ of (\ref{ds1}).
The y-component is determined by the equation
\begin{equation}
y(t)=h(t,s(t),\dot{s}(t)) \label{des0}%
\end{equation}

The pair of equations (\ref{de1}) and (\ref{des0}) constitute a folding of the
differential system (\ref{ds1}). Equation (\ref{des0}) is a passive equation
since it merely evaluates a function once $s(t)$ is known.

In general, the differential equation (\ref{de1}) may be no easier to deal
with than the system (\ref{ds1}). However, if (\ref{de1}) does not have
chaotic flows then, given the passive nature of (\ref{des0}) we may conclude
that (\ref{ds1}) also does not have chaotic flows. An example serves to
illustrate the considerable transformation that may occur as a result of
folding. Consider the system
\begin{equation}
\left\{
\begin{array}
[c]{l}%
\dot{x}=x^{2}+e^{a\sin bt}y+c\\
\dot{y}=-2e^{-a\sin bt}x^{3}-2xy-aby\cos bt
\end{array}
\right.  \label{des1}%
\end{equation}
where we may assume that the constants $a,b,c\not =0.$ This system is defined
in the entire plane. Here $f(t,x,y)=x^{2}+e^{a\sin bt}y+c$ has partial
derivatives%
\[
f_{t}=ab(\cos bt)e^{a\sin bt}y,\quad f_{x}=2x,\quad f_{y}=e^{a\sin bt}%
\]

Also we solve the first equation of (\ref{des1}) for $y$
\begin{equation}
y=e^{-a\sin bt}(\dot{x}-x^{2}-c)=h(t,x,\dot{x}). \label{des2}%
\end{equation}

Now (\ref{de2}) (or direct differentiation with respect to $t$ and
substitution) yields%
\[
\overset{..}{x}=ab(\dot{x}-x^{2}-c)\cos bt+\dot{x}(2x)-2x^{3}-2x(\dot{x}%
-x^{2}-c)-ab(\dot{x}-x^{2}-c)\cos bt=2cx
\]

Suppose that at $t=0$ the values $x(0)=x_{0}$ and $y(0)=y_{0}$ are given. If
$s(t)$ is a solution of the linear initial value problem%
\begin{equation}
\overset{..}{s}=2cs,\quad s(0)=x_{0},\ \dot{s}(0)=x_{0}^{2}+y_{0}+c
\label{nce}%
\end{equation}
then a flow $(x(t),y(t))$ of the system is obtained with $x(t)=s(t)$ and
$y(t)$ given passively by (\ref{des2}). Since the second-order differential
equation in (\ref{nce}) does not have chaotic flows, neither does the system
(\ref{des1}). In fact, the general solution of (\ref{des1}) can be found in
terms of elementary functions using (\ref{des2}) and (\ref{nce}).

To obtain systems whose folding will not be chaotic we solve a partial inverse
problem, starting with a desired function $\phi$ for (\ref{de1}) and a given
$f$ from (\ref{ds1}), then determine the suitable function $g.$

Assume that $f$ is semi-invertible with $h$ as in (\ref{des}). Then the right
hand side of (\ref{de2}) may be written as follows as a function of
$t,u=x,w=\dot{x}$%
\begin{equation}
f_{t}(t,u,h(t,u,w))+f(t,u,h(t,u,w))f_{u}(t,u,h(t,u,w))+g(t,u,h(t,u,w))f_{v}%
(t,u,h(t,u,w))=\phi(t,u,w) \label{ic1}%
\end{equation}

Now suppose, inversely, that the function $\phi(t,u,w)$ is prescribed and we
wish to determine a function $g$ such that (\ref{ic1}) is satisfied. Define
\begin{equation}
g(t,u,v)=\frac{1}{f_{v}(t,u,v)}\left[  \phi(t,u,f(t,u,v))-f(t,u,v)f_{u}%
(t,u,v)-f_{t}(t,u,v)\right]  \label{gc}%
\end{equation}
provided that $f_{v}(t,u,v)\not =0.$ It is evident that if the above $g$ is
inserted into (\ref{de2}) with $u=\dot{x}$ and $v=h(t,x,\dot{x})$ then
cancellations yield (\ref{de1}). As a check if $f(t,u,v)=u^{2}+e^{a\sin
bt}v+c$ as in (\ref{des1}) and $\phi(t,u,w)=2cu$ then (\ref{gc}) yields the
function $g$ in (\ref{des1}).

The preceding observations imply the following.

\begin{proposition}
\label{cds}Let $f$ be a semi-invertible function with a semi-inversion $h$ as
in (\ref{des}) and let $\phi$ be a function such that (\ref{de1}) has no
chaotic flows. If $f_{v}(t,u,v)\not =0$ for all $t,u,v$ and $g$ is given by
(\ref{gc}) then (\ref{ds1}) has no chaotic flows.
\end{proposition}

The following are straightforward corollaries.

\begin{corollary}
\label{c}Let $f$ be a semi-invertible function with a semi-inversion $h$ as in
(\ref{des}) and let $\phi(t,u,w)=\sigma(u,w)$ be independent of $t$ on an
interval $I.$ If $f_{v}(t,u,v)\not =0$ for all $t,u,v$ and%
\[
g(t,u,v)=\frac{1}{f_{v}(t,u,v)}\left[  \sigma(u,f(t,u,v))-f(t,u,v)f_{u}%
(t,u,v)-f_{t}(t,u,v)\right]
\]
then (\ref{ds1}) does not have chaotic flows.
\end{corollary}

\begin{corollary}
Let $f$ be a semi-invertible function with a semi-inversion $h$ as in
(\ref{des}) and let $\phi(t,u,w)=Au+Bw+C(t)$ where $A,B$ are constants, not
both zero, and $C(t)$ is a continuous function on an interval $I$. If
$f_{v}(t,u,v)\not =0$ for all $t,u,v$ and%
\[
g(t,u,v)=\frac{1}{f_{v}(t,u,v)}\left[  Au+(B-f_{u}(t,u,v))f(t,u,v)-f_{t}%
(t,u,v)+C(t)\right]
\]
then (\ref{ds1}) does not have chaotic flows.
\end{corollary}

We close by observing that folding may be useful in establishing the
occurrence of chaos also if it is known that the scalar equation (\ref{de1})
has chaotic flows. For example, the system%
\begin{equation}
\left\{
\begin{array}
[c]{l}%
\dot{x}=x+y+\delta(t)\\
\dot{y}=ax-bx^{3}-cy
\end{array}
\right.  \quad b,c>0\label{dufs}%
\end{equation}
folds to the scalar equation%
\begin{equation}
\overset{..}{s}+(c-1)\dot{s}-(a+c)s+bs^{3}=\delta^{\prime}(t)+c\delta
(t)\label{dufg}%
\end{equation}
together with the passive equation $y(t)=\dot{x}-x-\delta(t).$ If $\delta(t)$
is a solution of the linear differential equation
\begin{equation}
\delta^{\prime}(t)+c\delta(t)=A\sin\omega t\label{lin}%
\end{equation}
then (\ref{dufg}) is the periodically forced Duffing's equation which has
chaotic flows for a range of parameters \cite{ST}. The general solution of
(\ref{lin}) is%
\[
\delta(t)=Ke^{-ct}-\frac{A}{c^{2}+\omega^{2}}(c\sin\omega t+\omega\cos\omega
t)
\]
where $K$ is a constant of integration. For this $\delta$, (\ref{dufs}) is
chaotic for a range of parameter values. Of course, if $A=0$ then
$\delta(t)=Ke^{-ct}$ is a solution of the linear equation $\delta^{\prime
}(t)+c\delta(t)=0$. In this case, (\ref{dufg}) is autonomous so by Corollary
\ref{c} there are no chaotic flows.

\end{document}